\documentclass[12pt,reqno]{amsart}
\usepackage{amssymb,amsthm} 
\usepackage[utf8]{inputenc}

\advance\hoffset by -0.5truein
\let\<\langle
\let\>\rangle
\def\oo#1{\mathrel { {}_{(#1)}}}

\def\Cur{\mathop {\fam0 Cur}\nolimits}
\def\Hom{\mathop {\fam0 Hom}\nolimits}
\def\Cend{\mathop {\fam0 Cend}\nolimits}
\def\id{\mathop {\fam0 id}\nolimits}
\def\As{\mathrm{As}}
\def\Lie{\mathrm{Lie}}
\def\Conf {\mathop {\fam0 Conf }\nolimits}

\def\Com{\mathrm{Com}}

\newtheorem{theorem}{Theorem}

\newtheorem{proposition}{Proposition}
\newtheorem{corollary}{Corollary}

\theoremstyle{definition}
\newtheorem{definition}{Definition}
\newtheorem{remark}{Remark}
\newtheorem{example}{Example}

\title[Standard bases for the universal associative envelopes]{Standard bases for the universal associative conformal envelopes of Kac--Moody conformal algebras}

\thanks{The work is supported by Mathematical Center in Akademgorodok}

\author[P.S. Kolesnikov, R.A. Kozlov]{P.S. Kolesnikov$^{1)}$, R.A. Kozlov$^{1)2)}$}

\address{${}^{1)}$ Sobolev Institute of Mathematics, Novosibirsk, Russia}
\address{${}^{2)}$ Novosibirsk State University, Novosibirsk, Russia}

\begin{document}
\begin{abstract}
We study the universal enveloping associative conformal algebra 
for the central extension of a current Lie conformal algebra at the 
locality level $N=3$. A standard basis of defining relations for this algebra 
is explicitly calculated. 
As a corollary, we find a linear basis of the free commutative conformal algebra 
relative to the locality $N=3$ on the generators.
\end{abstract}


\maketitle

\section{Introduction}

Conformal algebras also known as Lie vertex algebras were introduced in 
\cite{K1} as an algebraic tool to study the singular part of the operator 
product expansion (OPE) of chiral fields in 2-dimensional conformal field theory 
coming back to \cite{BPZ}.
From the categorical point of view, a conformal algebra is just an algebra in 
the appropriate (pseudo-tensor) category 
$\mathcal M^*(\mathbb C[\partial ])$ of modules over the polynomial algebra 
$\mathbb C[\partial ]$ in one variable \cite{BDK}.
The pseudo-tensor structure (see \cite{BD}) reflects the main features of 
multi-linear maps in the category of linear spaces: 
composition, identity, symmetric structure. 
These features are enough to define the basic notions like what is an algebra 
(associative, commutative, Lie, etc.), homomorphism, ideal, representation, 
module, cohomology. 
Therefore, the notion of a conformal 
algebra is a natural expansion of the notion of an ``ordinary'' algebra over 
$\mathbb C$ to the pseudo-tensor category
$\mathcal M^*(\mathbb C[\partial ])$. 
Namely, as an ordinary algebra is a linear space equipped with 
a bilinear product, a conformal algebra is a 
$\mathbb C[\partial ]$-module $V$ equipped with a
$\mathbb C[\partial ]$-bilinear map (pseudo-product)
\[
*: V\otimes V \to \mathbb C[\partial ]^{\otimes 2}\otimes_{\mathbb C[\partial ]} 
V.
\]
A more convenient presentation for the operation $*$
uses the language of a $\lambda $-product or a family of $n$-products
for all integer $n\ge 0$ (\cite{K1}, see also Section~\ref{sec:Conf}).

Conformal algebras representing the singular part of OPE in vertex algebras are 
Lie algebras in the category $\mathcal M^*(\mathbb C[\partial ])$, i.e., Lie 
conformal algebras. 
For example, if $\mathfrak g$ is a Lie algebra then the 
free module $\mathbb C[\partial ]\otimes \mathfrak g$ equipped 
with the pseudo-product 
$a*b = (1\otimes 1)\otimes _{\mathbb C[\partial ]} [a,b]$, 
$a,b\in \mathfrak g$, is a Lie conformal algebra denoted 
$\Cur \mathfrak g$ (current conformal algebra). 
If $\<\cdot|\cdot \>$ is a bilinear symmetric invariant form on 
$\mathfrak g$ then $\Cur\mathfrak g $  has a 1-dimensional 
central extension $K(\mathfrak g)$ defined by
\[
a*b = (1\otimes 1)\otimes _{\mathbb C[\partial ]} [a,b]
- (\partial \otimes 1)\otimes _{\mathbb C[\partial ]} \<a|b\>e, 
\]
where $e$ is a central element and $\partial e=0$.
For example, in the Kac--Moody vertex algebra $V(\mathfrak g)$ \cite{FBZ} the 
singular part of the OPE on the generating 
fields is described by this particular structure $K(\mathfrak g)$
called a Kac--Moody conformal algebra. 

As in the case of ordinary algebras, an associative conformal algebra $C$ turns into 
a Lie one with respect to the commutator
$[a*b] = (a*b)-(\tau \otimes _{\mathbb C[\partial ]} 1)(b*a)$, $a,b\in C$, 
where $\tau $ is the switching map on $\mathbb C[\partial ]^{\otimes 2}$. 
However, not all Lie conformal algebras embeds into associative ones in this way 
\cite{Ro2000}. This is an open problem whether every finite (i.e., finitely 
generated as a $\mathbb C[\partial]$-module) 
Lie conformal algebra embeds into an associative conformal algebra 
with respect to conformal commutator. Even for the class of quadratic conformal 
algebras  \cite{Xu} (see also \cite{HongWu17}) 
it remains unknown in general if every such Lie conformal algebra embeds into an 
appropriate associative one. 

A routine way to solve this kind of problems is to construct a universal 
envelope. 
In general, such an algebra is defined by generators and relations.
For  a Lie conformal algebra, there exists a lattice of universal enveloping 
associative conformal algebras, 
each related to an (associative) locality bound on the generators
(\cite{Ro2000}, see also Section~\ref{subsec_env}).
In order to prove (or disprove) the embedding of a Lie conformal algebra into 
its universal enveloping associative conformal algebra 
one needs to know the normal form of elements in the last algebra.

A general and powerful method for finding normal forms in an algebra defined by 
generators and relations is to calculate a standard (or Gr\"obner--Shirshov) 
basis of defining relations. The idea goes back to  Newmann's Diamond Lemma 
\cite{New53}, see also 
\cite{Bok72, Berg7x}. 
In the recent years, the Gr\"obner--Shirshov bases theory was developed to serve 
the problem of combinatorial analysis of various algebraic structures, see 
\cite{BokChen14}. 
For associative conformal algebras it was initially invented in \cite{BFK2000}, 
later developed in \cite{NC2017} and \cite{Ko2020}. 
In this paper, we use the last approach exposed in a form convenient for actual 
computation: we consider defining relations in a conformal algebra as rewriting 
rules on a module over an appropriate associative algebra (the 
Gr\"obner--Shirshov basis of the last algebra is known).

A series of particular observations made in \cite{KolIJAC}, \cite{KKP} 
shows that for all considered examples of quadratic Lie conformal algebras $L$
it is enough to consider universal associative conformal envelopes $U$ relative 
to the locality bound $N=3$ to get an injective mapping $L\curvearrowright U$.
This is one of the reasons why we focus on the locality bound $N=3$
for the envelopes of current conformal algebras as they are particular examples 
of quadratic conformal algebras.

The main purpose of this paper is to find a standard (Gr\"obner--Shirshov) basis 
of defining relations 
for the universal enveloping associative conformal algebra 
of a Kac--Moody conformal algebra at locality level $N=3$.
As a corollary, we get an analogue of the Poincar\'e--Birkhoff--Witt Theorem 
(PBW-Theorem) stating that the associated graded 
conformal algebra obtained from the universal envelope of a current Lie 
conformal algebra with respect to 
the natural filtration is isomorphic to the free commutative conformal algebra. 
Note that the classical PBW-Theorem may be interpreted as a conformal one at the 
locality level $N=1$: for a Lie algebra $\mathfrak g$, its ``ordinary'' 
universal envelope $U(\mathfrak g)$ gives rise to the conformal algebra $\Cur 
U(\mathfrak g)_0$ which is exactly the universal enveloping conformal algebra of 
$\Cur \mathfrak g$ with $N=1$ (here $U(\mathfrak g)_0$ is the augmentation ideal 
of $U(\mathfrak g)$).

There are several reasons for studying the universal envelopes of 
$\Cur \mathfrak g$ at higher locality than $N=1$. 

First, the $N=1$ envelope $\Cur U(\mathfrak g)$ does not reflect the homological 
properties of $\Cur \mathfrak g$. For example, if $\mathfrak g$ is a simple 
finite-dimensional Lie algebra then 
the second cohomology group
$H^2(\Cur\mathfrak g, \Bbbk )$ is one-dimensional \cite{BKV}.
The corresponding central extension is the Kac--Moody conformal algebra 
$K(\mathfrak g)$ representing the singular part of the Kac--Moody vertex algebra 
\cite{FBZ}.
On the other hand, it is easy to find that the second Hochschild cohomology 
group of $\Cur U(\mathfrak g)_0$ with coefficients in the trivial 1-dimensional module 
is zero: there are no nontrivial central extensions. Our results  show that 
the universal enveloping associative conformal algebras for $\Cur \mathfrak g$ 
at locality level $N=2,3$ do have a nontrivial central extension which is 
exactly the universal envelope of $K(\mathfrak g)$.

The second reason is related with Poisson algebras. Assume $P$ is an ordinary 
commutative algebra with a Poisson bracket 
$\{\cdot ,\cdot \}$. Then $\Cur P$ may be considered as a Lie conformal algebra 
since $P$ is a Lie algebra relative to the Poisson bracket. There is a conformal 
representation of $\Cur P$ on itself given by the rule 
\[
(a\oo\lambda f) = \{a,f\} +\lambda af,\quad a,f\in P.
\]
The study of this representation provides a way to get new results 
in (quadratic) conformal algebras as well as in Poisson algebras 
\cite{KolIJAC, KKP}.
The conformal linear operators $\rho(a)\in \Cend(\Cur P)$,
$\rho(a)_\lambda : f\mapsto (a\oo\lambda f)$, 
are local to each other, and the locality bound is $N=3$.
Indeed, according to the definition of a conformal representation
\cite{CK, K1} we have 
\begin{multline*}
(\rho(a)\oo\lambda \rho(b))_\mu f
= a\oo\lambda (b\oo {\mu-\lambda} f) \\
= \{a,\{b,f\}\} +\lambda a\{b,f\} 
+ (\mu-\lambda )\{a,bf\} +\lambda (\mu-\lambda )abf,
\end{multline*}
for $a,b,f\in P$. If $abf\ne 0$ then right-hand side is of degree 2 in $\lambda 
$ that means $N(\rho(a),\rho(b))=3$ in $\Cend (\Cur P)$.
Therefore, the corresponding associative envelope belongs to the class of 
envelopes with locality $N=3$.

The third reason to study the case $N=3$ comes from the following relation 
between commutative conformal and Novikov algebras. 
Suppose $C$ is a commutative conformal algebra and $M$ is a subset of $C$ such 
that $N_C(a,b)\le 3$ for all $a,b\in M$. 
Then $M$ generates an ordinary (nonassociative) subalgebra $N(M)$ in the space 
$C$ considered relative to the single product $x\circ y = x\oo{1} y$. Indeed, 
all elements of $N(M)$ are local to each other with locality bound $3$. 
Moreover, the following relations hold:
\[
\begin{gathered}
 (x\circ y)\circ z -x\circ (y\circ z) =  
 (x\circ z)\circ y -x\circ (z\circ y), \\
 x\circ (y\circ z) = y\circ (x\circ z),
\end{gathered}
\]
for all $x,y,z\in N(M)$.
These identities are known to define the variety of Novikov algebras initially 
appeared in  \cite{BN89}, \cite{GD89}. 
In order to perform a systematic study of this relation, one needs to know the 
structure of the universal object in the category of commutative conformal 
algebras with locality bound $N=3$ on the generators. 

For all these reasons, we study the universal enveloping associative conformal 
algebras for Kac--Moody conformal algebras $K(\mathfrak g)$ relative to the 
locality level $N=3$. The corollaries of the main result of the paper 
(Theorem~\ref{thm:GSB}) allow us to get the structure of the universal envelopes 
for current Lie conformal conformal algebras at $N=3$ and also describe the free 
commutative conformal algebra at the same locality level. Practically, we find a 
standard (Gr\"obner--Shirshov) basis of defining relations for these conformal 
algebras and derive an analogue of the PBW-Theorem.

\section{Preliminaries in conformal algebras}\label{sec:Conf}

The definition of a conformal algebra as an algebra in an appropriate pseudo-tensor 
category \cite{BDK}
corresponds to the convenient algebraic approach using $\lambda $-brackets \cite{K1}
if it is presented in terms of operads associated with linear algebraic groups \cite{Ko2009}.

Let $G$ be a linear algebraic group over a field $\Bbbk $ of characteristic 
zero, 
and let $H_G =\Bbbk [G]$ be the Hopf algebra 
of regular functions on~$G$. 
For every $H_G$-module $V$ there is a non-symmetric operad (let us denote it 
$V_G$) 
defined as follows. Given $n\in \{1,2,\dots \}$, set 
\[
 V_G(n) = \{ f: G^{n-1}\to \Hom (V^{\otimes n}, V) \mid \text{$f$ is regular and 
3/2-linear} \}
\]
The condition of regularity means that $f$ may be presented by a polynomial 
function with coefficients 
in the space $\Hom (V^{\otimes n}, V)$ of $\Bbbk $-polylinear maps on $V$, and 
the 3/2-linearity (sesqui-linearity)
may be expressed as 
\[
 f(\lambda_1,\dots , \lambda_{n-1}): (v_1, \dots, h(x)v_i, \dots , v_n) 
\mapsto 
 \begin{cases}
  h(\lambda_i^{-1})v, & i=1,\dots, n-1 , \\
  h(\lambda_{n-1}\dots \lambda_1x) v, & i=n,
\end{cases}
\]
for 
$v=f(\lambda_1,\dots , \lambda_{n-1})(v_1,\dots, v_n)$,
$\lambda_i\in G$, $v_i\in V$, 
$h(x)\in H_G$ (here $x$ is a variable ranging in $G$). 
In particular, $V_G(1)$ is the space of all $H_G$-linear transformations of $V$
thus contains the identity map $\id$.

The composition rule in $V_G$ is defined by the following partial composition. 
If $f\in V_G(n)$, $g\in V_G(m)$, $i\in \{1,\dots , n\}$ then 
\[
 f\circ_i g = f(\id,\dots , \underset{i}{g}, \dots , \id) \in V_G(n+m-1)
\]
acts as
\begin{multline}\label{eq:Gmod_comp}
f\circ_i g: (\lambda_1,\dots, \lambda_{i-1},\mu_1,\dots , \mu_{m-1}, 
\lambda_i,\dots , \lambda_{n-1}) \\
\mapsto 
 f(\lambda_1,\dots,\lambda_{i-1}, \lambda_i\mu_{m-1}\dots \mu_{1}, 
\lambda_{i+1},\dots,\lambda_{n-1})\circ_i
 g(\mu_1,\dots , \mu_{m-1}),
\end{multline}
for $\lambda_i,\mu_j\in G$,
where $\circ_i$ in the right-hand side stands for the ordinary partial 
composition 
of polylinear maps.
In particular, for $i=n$ the partial composition is equal to
\[
f(\lambda_1,\dots, \lambda_{n-1})\circ_n g(\mu_1,\dots , \mu_{m-1}).
\]
It is easy to see that the resulting maps are indeed regular and 3/2-linear.

One may easily check that the partial composition in $V_G$ defined above meets 
the sequental, parallel, and unit axioms \cite[Definition 3.2.2.3]{BD_AlgOP}
and thus this is indeed an non-symmetric operad with a well-defined 
composition rule 
\[
 \gamma^{r}_{m_1,\dots , m_r}: V_G(r)\otimes V_G(m_1)\otimes \dots \otimes 
V_G(m_r)
 \to V_G(m_1+\dots + m_r).
\]

Suppose the group $G$ is abelian. Then $V_G(n)$ has a natural action of the 
symmetric group $S_n$ defined in the following way. 
 If $f\in V_G(n)$ and $(1i)$, $i=2,\dots, n$, is a transposition in $S_n$ then 
\[
 f^{(1i)}(\lambda_1,\dots, \lambda_{n-1}) = f(\lambda_i,\lambda_2,\dots , 
\underset{i}{ \lambda_1}, \dots \lambda_{n-1})^{(1i)}
\]
for $i<n$ (here the action of $(1i)$ in the right-hand side is just the 
permutation of arguments in a polylinear map). 
For $i=n$, the definition is slightly more complicated: 
if $f$ is presented by a polynomial function 
\[
 f = \sum\limits_i f_i(x_1,\dots, x_{n-1})\varphi_i,\quad f_i\in H_G^{\otimes 
n-1}\simeq H_{G^{n-1}}, \ \varphi_i \in \Hom (V^{\otimes n}, V)
\]
then 
$f^{(1n)} (\lambda_1,\dots, \lambda_{n-1})$ is given by 
\[
 \sum\limits_i f'_i(x)\varphi_i^{(1n)},
\]
where the each $f'_i(x)\in H_G$ is the regular function  
$f_i((\lambda_1\dots \lambda_{n-1}x)^{-1}, \lambda_2,\dots, \lambda_{n-1})$.

The composition rule $ \gamma^{r}_{m_1,\dots , m_r}$  is equivariant (see, e.g., 
\cite[Definition 5.2.1.1]{BD_AlgOP})
since the structure obtained is equivalent to the structure of an $H_G$-module 
operad defined over 
a cocommutative Hopf algebra \cite{BDK}. Namely, one may identify a map 
\[
 F: V^{\otimes n} \to H_G^{\otimes n}\otimes _{H_G} V\simeq
 H_G^{\otimes n-1}\otimes V\simeq
 H_{G^{n-1}}\otimes V,
\]
with 
\[
 f(\lambda_1,\dots, \lambda_{n-1}) = F(\lambda_1^{-1},\dots, \lambda_n^{-1}).
\]

Recall that if $\mathcal O$ is a (symmetric) operad then a morphism $\mathcal 
O\to V_G$ defines an algebraic structure on 
the space $V$.
In the trivial case $G=\{e\}$, $H_G=\Bbbk $, the operad $V_G$ coincides with the 
operad of polylinear maps on the linear space $V$, and thus 
a morphism $\mathcal O\to V_G$ defines the ordinary notion of an $\mathcal 
O$-algebra over $\Bbbk $, a space equipped with a family of polylinear 
operations.

The classes of conformal \cite{K1} and $\mathbb Z$-conformal \cite{GKK} algebras  naturally 
appear in the next step, if we choose $G$ to be 
a connected linear algebraic group of dimension 1 
(the affine line and $GL_1$, respectively). 
For a non-connected group $G$ with the identity component denoted $G^0$, the 
structure of a conformal
algebra over $G$ is naturally interpreted as 
a $G/G^0$-graded conformal algebra over $G^0$ \cite{Kol2013_art}.
If $G=\mathbb A_1 = (\Bbbk, +)$, $H_G=\Bbbk [\partial ]$ (the variable is 
traditionally denoted by $\partial $), 
then a morphism $\mathcal O\to V_G$ defines a $\mathcal O$-conformal algebra 
structure on 
a $\Bbbk [\partial ]$-module~$V$.

For example, 
if $\mathcal O=\As$ is the operad governing the variety of associative algebras 
(generated by $\mu = x_1x_2\in \As(2)$
modulo the relation $\mu\circ_1\mu = \mu \circ_2 \mu $)
then an associative conformal algebra structure on a 
$\Bbbk [\partial]$-module $V$ 
is given by an image of $\mu $, a map $f = (\cdot\oo\lambda \cdot): V\otimes V 
\to \Bbbk[\lambda ]\otimes V$
which is 3/2-linear 
\[
 (\partial v\oo\lambda u) = -\lambda (v\oo\lambda u), \quad 
 (v\oo\lambda \partial u) = (\partial +\lambda) (v\oo\lambda u),
\]
for $u,v\in V$,
and associative in the sense that 
\[
 (f\circ_2 f)(\lambda,\mu) = (f\circ_1 f)(\lambda, \mu).
\]
By \eqref{eq:Gmod_comp}, the latter means
\begin{equation}\label{eq:Conf_ass}
 (u\oo{\lambda } (v\oo\mu w)) = ((u\oo\lambda v)\oo{\lambda+\mu } w)
\end{equation}
(to compute the right-hand side, put $\mu_1=\lambda $ and $\lambda_1=\mu$ in 
\eqref{eq:Gmod_comp}).

According to the same scheme,
a Lie conformal algebra structure on a $\Bbbk [\partial ]$-module $V$
is a morphism from the operad $\Lie $ governing the variety of Lie algebras to 
$V_{\mathbb A_1}$.
To define such a morphism, it is enough to fix a 3/2-linear map 
$\mu \in V_{\mathbb A_1}(2)$
\[
 \mu = {[\cdot\oo\lambda \cdot ]}: V\otimes V\to \Bbbk [\lambda ]\otimes V
\]
such that 
$\mu^{(12)}=-\mu $ and 
$(\mu \circ_ 2 \mu ) - (\mu \circ_2 \mu )^{(12)}
 = (\mu \circ_1 \mu )$. 
The last two relations represent anti-commutativity and Jacobi identity, 
respectively:
\[
 [u\oo{-\partial -\lambda} v ] = - [v\oo\lambda u ],
\]
\[
 [u\oo\lambda [v\oo\mu w]] - [v\oo\mu [u\oo\lambda w]] = [ [u\oo\lambda v] 
\oo{\lambda+\mu} w ],
\]
$u,v,w\in V$.

In the sequel, we will use the notation $(\cdot \oo\lambda \cdot)$ for 
the operation on an associative conformal algebra and 
$[\cdot \oo\lambda \cdot ]$ for Lie conformal algebras. 

Since there is a morphism of operads $(-):\Lie \to \As$ sending $\mu $ to $f - 
f^{(12)}$, every associative conformal algebra turns into 
a Lie conformal algebra relative to the operation 
\[
[u\oo\lambda v ] = (u\oo\lambda v ) - (v\oo{-\partial -\lambda } u).
\]
For an associative conformal algebra $V$ defined via a morphism of operads 
$\As\to V_{\mathbb A_1}$, let $V^{(-)}$ stand 
for the Lie conformal algebra obtained as a composition 
$\Lie \overset{(-)}{\to} \As \to V_{\mathbb A_1}$.

The property of a commutator to be a derivation on an associative algebra may 
also be expressed as a relation in $\As(3)$. 
Being translated to conformal algebras it turns into the following identity on 
an associative conformal algebra $V$:
\begin{equation}\label{eq:Conf_JacComm}
(u\oo\lambda (v\oo\mu w))-(v\oo\mu (u\oo\lambda w))
= ([u\oo\lambda v]\oo{\lambda+\mu} w),
\quad u,v,w\in V.
\end{equation}

As in the case of ordinary algebras, $V\mapsto V^{(-)}$ is a functor from the 
category of associative conformal algebras 
to the category of Lie algebras. In contrast to the case of ordinary algebras, 
this functor does not have a left adjoint one when considered on the 
entire category of associative conformal algebras.
However, if we restrict the class of associative conformal algebras by means of 
locality on the generators (\cite{Ro2000}, see Section~\ref{subsec_env} for details) then there 
is an analogue 
of the universal enveloping associative algebra for Lie conformal 
algebras.

In terms of ``ordinary'' algebraic operations, a conformal algebra 
is a linear space $V$ equipped with a linear operator $\partial $, the generator 
of $H_{\mathbb A_1}=\Bbbk [\partial ]$, 
and a series of bilinear operations $(\cdot \oo n\cdot )$, 
$n\in \mathbb Z_+$,
given by 
\[
(u\oo\lambda v) = \sum\limits _{n\ge 0} \frac{\lambda^n}{n!} (u\oo n v),\quad 
u,v\in V.
\]
These operations are called $n$-products. They have to satisfy the following 
properties:
\begin{itemize}
    \item [(C1)] For every $u,v\in V$ there exists $N=N(u,v)$ such that $(u\oo n 
v)=0$ for all $n\ge N$; 
    \item [(C2)] $(\partial u\oo n v) = -n (u\oo {n-1} v)$; 
    \item [(C3)] $(u\oo n\partial v) = \partial (u\oo n v) + n(u\oo{n-1} v )$.
\end{itemize}
The property (C1) is known as the locality axiom, (C2) and (C3) represent 
3/2-linearity. 
For every conformal algebra $V$, the locality function $N_V$ 
is a map $V\times V\to \mathbb Z_+$ such that
$u\oo{n}v=0$ for every $u,v\in V$ and $n\ge N_V(u,v)$.

A conformal algebra $V$ is associative if 
\[
( u\oo n( v\oo m w)) = \sum\limits_{s=0}^n \binom{n}{s} 
((u\oo{n-s} v)\oo{m+s} w)
\]
for all $u,v,w\in V$ and $n,m\in \mathbb Z_+$. 
In a similar way, one may rewrite the identities defining 
the class of Lie conformal algebras.

Given a set $X$ and a function $N:X\times X\to \mathbb Z_+$, 
there exists a unique (up to isomorphism) associative conformal algebra denoted 
$\Conf(X,N)$ which is universal among all associative 
conformal algebras $V$ generated by $X$ such that $N_V(x,y)\le N(x,y)$ for all 
$x,y\in X$ \cite{Ro1999}.
The details of the construction of $\Conf(X,N)$ are stated in 
Section~\ref{subsec_free}.

\section{Gr\"obner--Shirshov bases for associative conformal algebras}

\subsection{Rewriting system and standard bases for associative algebras}
In this section, we briefly describe the well-known technique of standard bases 
(Gr\"obner--Shirshov bases) 
in associative algebras in order to fix the notations.
The usual exposition of this technique requires a proper ordering of the 
monomials. However, the core statements laying in the foundation of the approach 
do not need a monomial ordering. 

Let $B$ be a set and let $B^*$ stand for the set of all words in $B$ (including 
the empty word). 
The free associative algebra (with a unit) over the field $\Bbbk $ generated by 
$B$ is denoted $\Bbbk\<B\>$.
Suppose $\Sigma $ is a family of pairs $(u,f)$ called rewriting rules, where 
$u\in B^*$, $f\in \Bbbk \<B\>$. We will write a pair like this 
as $(u\to f)$ since the family $\Sigma $ determines an oriented graph $\mathcal 
G(B,\Sigma) $ as follows. 
The vertices of $\mathcal G(B,\Sigma )$ are the elements of $\Bbbk \<B\>$; two 
vertices $g$ and $h$ are connected with an edge 
($g\to h$)
if and only if there is a rewriting rule $u\to f$ in $\Sigma $ and a summand of 
the form $\alpha w$ in $g$ ($\alpha \in \Bbbk^\times$, $w\in B^*$)
such that 
\[
 w = v_1uv_2, \quad h=g-\alpha v_1(u-f)v_2
\]
for some $v_1,v_2\in B^*$.
In other words, $h$ is obtained from $g$ by replacing an occurrence of 
the subword $u$ with the polynomial~$f$.

The graph $\mathcal G(B,\Sigma )$ splits into connected components (in the 
non-oriented sense) which explicitly correspond to the elements 
of the quotient $\Bbbk \<B\mid \Sigma \>=\Bbbk \<B\>/(\Sigma )$, 
where $(\Sigma )$ stands for the 
ideal in $\Bbbk \<B\>$ generated by all $u-f$ for $(u\to f)\in \Sigma $.
In some cases, there is a way to check algorithmically whether two vertices 
$g,h\in \Bbbk\<B\>$ belong to the same connected component of 
$\mathcal G(B,\Sigma )$, i.e., if the images of $f$ and $g$ are equal in $\Bbbk 
\<B\mid \Sigma \>$.

An oriented graph is called a rewriting system if there are no infinite oriented 
paths (in particular, no oriented cycles). 
In a rewriting system, for every vertex $g$ there is a nonempty set $T(g)$ of 
terminal vertices $t$ attached to $g$, i.e., 
such that there is a path $g\to \dots \to t$, but there are no edges originated at~$t$.
A rewriting system is confluent if for every vertex $g$ the set $T(g)$ contains 
a single vertex.

\begin{definition}
 A family of rewriting rules $\Sigma $ in the free associative algebra $\Bbbk 
\<B\>$ 
 is a standard basis (Gr\"obner--Shirshov basis, GSB) if $\mathcal G(B,\Sigma 
)$ 
 is a confluent rewriting system.
\end{definition}

Obviously, if $\Sigma $ is a GSB then every connected component of $\mathcal 
G(B,\Sigma )$ 
has a unique terminal vertex which is a linear combination of terminal (reduced) 
words.
This combination is called a normal form of an element in $\Bbbk \<B\mid \Sigma \>$: 
two polynomials 
$g$ and $h$ in $\Bbbk \<B\>$ represent the same element of $\Bbbk \<B\mid \Sigma \>$ 
if and only if 
their normal forms coincide. Therefore, the images of terminal words form a 
linear basis 
of $\Bbbk \<B\mid \Sigma \>$. 

The most natural way to guarantee that $\mathcal G(B,\Sigma )$ is a rewriting 
system 
is to make the set $B^*$ well-ordered relative to an order $\le $ such that 
$u\le v$ implies $wu\le wv$ and $uw\le vw$ for all $u,v,w\in B^*$ 
(i.e., $\le $ is a monomial order), and $u>f$ for all $(u\to f)\in \Sigma $
(i.e., $u$ is greater than every monomial in $f$). 

To check the confluence of a rewriting system $\mathcal G(B,\Sigma )$ 
one may apply the Diamond Lemma originated to \cite{New53}. 
The latter states that a rewriting system is confluent if an only if 
for every ``fork'' (a pair of edges $w\to g_1$, $w\to g_2$) there exist a vertex 
$h$ 
and two oriented paths $g_1\to \dots \to h$, $g_2\to \dots \to h$.
If the rewriting system is $\mathcal G(B,\Sigma)$ then it is enough to check the 
Diamond condition for the following two kinds of forks: 
\begin{enumerate}
 \item For $u_1\to f_1$, $u_2\to f_2$ in $\Sigma $, 
 $u_1=v_1u_2v_2$, consider $w=u_1$, $g_1=f_1$, and $g_2= v_1f_2v_2$; 
 \item For $u_1\to f_1$, $u_2\to f_2$ in $\Sigma $,
 $u_1=v_1v$, $u_2=vv_2$, $v$ is a nonempty word, consider 
 $w=v_1vv_2$, $g_1=f_1v_2$, $g_2=v_1f_2$.
\end{enumerate}
In both cases, if there exit oriented paths $g_1\to \dots \to h$ and 
$g_2\to \dots \to h$ for an appropriate polynomial $h$ then we say that 
the composition of $u_1\to f_1$ and $u_2\to f_2$ relative to the word $w$ is 
confluent modulo $\Sigma $. 
Denote the polynomial $g_1-g_2$ by $(u_1\to f_1, u_2\to f_2)_w$.

\begin{theorem}[\cite{Berg7x, Bok72}]
Suppose a set of rewriting rules $\Sigma $ in
the free associative algebra $\Bbbk \<B\>$
defines a rewriting system $\mathcal G(B,\Sigma )$. 
If every composition of rewriting rules from $\Sigma $ is confluent modulo 
$\Sigma $ 
then $\mathcal G(B,\Sigma )$ is a confluent rewriting system, i.e., $\Sigma $ is 
a GSB.
\end{theorem}

Let $\Sigma $ respect a monomial order $\le $ on $B^*$. Then $\mathcal 
G(B,\Sigma )$ 
is a rewriting system and the confluence of a composition may be replaced with 
a more convenient condition.

\begin{corollary}[\cite{Bok72}]
If for every rewriting rules $u_1\to f_1$, $u_2\to f_2$ in $\Sigma $ having a 
composition relative to a word $w$
the polynomial $(u_1\to f_1, u_2\to f_2)_w$ may be presented as 
\begin{equation}\label{eq:CompTrivialOrder}
 \sum\limits_{i} \alpha_i w_i(u^{(i)}-f^{(i)}) w'_i, \quad \alpha_i\in \Bbbk , \ 
w_i,w_i'\in B^*,
\end{equation}
where 
$u^{(i)}\to f^{(i)}$ in $\Sigma $ and $w_iu^{(i)}w'_i < w$, 
then $\Sigma $ is a GSB.
\end{corollary}

In the actual computation, we will often apply the following trick to show the 
confluence of a fork 
$w\to g_1$, $w\to g_2$: find some paths $g_1\to \dots \to h_1$  and $g_2\to 
\dots \to h_2$ and then 
present $h_1-h_2$ in the form \eqref{eq:CompTrivialOrder}.

\subsection{Rewriting system for bimodules over associative algebras}
Let $A$ be an associative algebra (with a unit) and let $M$ be a bimodule over 
$A$. 
Suppose $A$ is generated by a subset $B\subset A$ as an algebra and $M$ 
is generated by a subset $Y$ as an $A$-module. 
Then $A$ is isomorphic to a quotient of the free associative algebra 
$\Bbbk\<B\>$ 
modulo an ideal generated by a set of defining relations $R\subset \Bbbk\<B\>$, 
i.e., 
\[
 A \simeq  \Bbbk\<B\mid R\>.
\]
Similarly, $M$ is a quotient of the free $A$-module  $A\otimes \Bbbk Y\otimes A$ 
generated by $Y$
modulo a family of defining relations $S$.
One may identify an element of $S$ with a noncommutative polynomial in the 
variables $B\cup Y$ 
which is linear in~$Y$.

The split null extension $A\oplus M$ is an associative algebra isomorphic to the 
quotient of the free 
algebra generated by $B\cup Y$ modulo the ideal generated by the union of $R$, 
$S$, and 
\[
yb_1\dots b_n z,\quad y,z\in Y,\ b_i\in B, \ n\ge 0.
\]
These relations reflect the properties of multiplication in $A\oplus M$: 
$M^2=0$.

\begin{remark}
 To consider left modules, it is enough to add relations $yb$, $y\in Y$, $b\in 
B$ 
 to reflect $MA=0$. 
\end{remark}

Suppose we may choose a monomial $u$ in each defining relation $u-f$ of $A\oplus 
M$ (up to a scalar multiple)
in such a way that the family $\Sigma $ of all rewriting rules $u\to f$ defines 
a rewriting system 
$\mathcal G(B\cup Y, \Sigma )$. 
Note that the defining relations of $A\oplus M$ are homogeneous relative to $Y$. 
All monomials that are of degree 
$\ge 2$ in $Y$ belong to the same connected component as zero, so it is enough 
to consider only the relations of degree 0 and 1 in $Y$, 
these are exactly the defining relations of $A$ and of $M$, respectively. 
Therefore, the confluence test needs to be applied 
to the forks started at a word $w$ which either belongs to $B^*$ or contains only one 
letter from~$Y$.
Hence, the compositions emerging in this rewriting system are exactly those 
described 
in \cite{KangLee}.

\subsection{Free associative conformal algebras}\label{subsec_free}

Recall the construction of a free associative conformal algebra $\Conf(X,N)$ 
generated by a set~$X$ 
relative to a given locality function $N: X\times X\to \mathbb Z_+$. 
From now on, denote by $H$ the polynomial algebra $\Bbbk [\partial ]$.

By definition, $\Conf (X,N)$ is an associative conformal algebra generated by 
$X$ which is 
universal in the class of all associative conformal algebras $C$ generated by 
$X$ such that 
the mutual locality of elements from $X$ in $C$ is bounded by~$N$. 
Namely, for every associative conformal algebra $C$ and for every 
map $\alpha: X\to C$ such that $N_C(\alpha (x),\alpha (y))\le N(x,y)$ for all $x,y\in 
X$
there exists unique homomorphism of conformal algebras 
$\varphi: \Conf(X,N)\to C$ 
such that $\varphi(x)=\alpha (x)$ for all $x\in X$.

\begin{proposition}[\cite{Ro1999}]\label{prop:FreeBasis}
The free associative conformal algebra $\Conf (X,N)$ 
is a free $H$-module with a basis 
\[
\begin{gathered}
 a_1\oo{n_1}(a_2\oo{n_2} (a_3\oo{n_3} \dots \oo{n_{k-1}}(a_k\oo{n_k} 
a_{k+1})\dots )), \\
  a_i\in X,\ 0\le n_i\le N(a_i,a_{i+1}),\ k\in \mathbb Z_+.
\end{gathered}
\]
\end{proposition}

\begin{remark}
In a similar way, one may define the free associative commutative conformal 
algebra $\Com\Conf(X,N)$
generated by a set $X$ relative to a locality function $N$ \cite{Ro2000}. 
However, there was no explicit description of a linear basis  of 
$\Com\Conf(X,N)$ for $N>1$. We will obtain such a description for $N=2,3$ as a 
byproduct in Section~\ref{sec:PBW-KM}.
\end{remark}

The conformal algebra $\Conf(X,N)$ may be presented in a more convenient form as a (left) 
module over an appropriate associative algebra \cite{Ko2020}.
Given a set $X$, let $A(X)$ denote the associative algebra generated by the set 
\[
 B=\{\partial \}\cup \{L_n^a,R_n^a \mid a\in X,\,n\in \mathbb Z_+\}
\]
relative to the defining relations
\begin{align}\label{eq:rel_A(X)}
 & L_n^a\partial - \partial L_n^a -n L_{n-1}^a, \\
 & R_n^a\partial - \partial R_n^a -n R_{n-1}^a, \\
 & R_n^aL_m^b - L_m^b R_n^a,
\end{align}
where  $a,b\in X$, $n,m\in \mathbb Z_+$.

The free associative conformal algebra $\Conf (X,N)$ is a left module over 
$A(X)$ if we define the action as follows:
\[
 L_n^a u = a\oo{n}u,\quad R_n^a u = \{u\oo{n} a\}
\]
for $a\in X$, $n\in \mathbb Z_+$, $u\in \Conf(X,N)$. 
Therefore, $\Conf (X,N)$ considered as a left $A(X)$-module is a homomorphic 
image 
the free left $A(X)$-module $M(X)$ generated by the set~$X$.
It is not hard to find explicitly the kernel of that homomorphism $M(X)\to \Conf 
(X,N)$.
 
Fix a function $N:X\times X\to \mathbb Z_+$ and 
consider the quotient $M(X,N)$ of $M(X)$ relative to the $A(X)$-submodule 
generated by the following elements:
\begin{align}\label{eq:mod_A(X)}
  & L_n^ab, \quad n\ge N(a,b),\\
  & R_m^ab - \sum\limits_{s= 
0}^{N(b,a)-m}(-1)^{m+s}\frac{1}{s!}\partial^{s}L_{m+s}^b a,
  \quad m\in \mathbb Z_+,
\end{align}
where $a,b\in X$.
Obviously, there is a homomorphism $M(X,N)\to \Conf (X,N)$ of $A(X)$-modules 
extending $x\mapsto x$.
This homomorphism is actually an isomorphism since 
 \eqref{eq:rel_A(X)} and \eqref{eq:mod_A(X)} imply the following relations in 
$M(X,N)$:
\begin{equation}\label{eq:mod_A(X)_locality}
L_n^aL_m^bu +\sum\limits_{q\ge 1} (-1)^q \binom{n}{q} L_{n-q}^aL_{m+q}^b u=0,
\end{equation}
where $a,b\in X$, $n\ge N(a,b)$, $m\in \mathbb Z_+$, $u\in M(X)$.

Consider the relations \eqref{eq:rel_A(X)}--\eqref{eq:mod_A(X)_locality}
as rewriting rules in such a way that the first monomial is always a principal 
one.
The terminal words in $M(X)$ of the rewriting system obtained are 
\[
 \partial ^s L_{n_1}^{a_1}L_{n_2}^{a_2}\dots L_{n_k}^{a_k}a_{k+1},\quad k\in 
\mathbb Z_+,\ a_i\in X, \ 0\le n_i<N(a_i,a_{i+1}),\ s\in \mathbb Z_+. 
\]
The images of these words in $\Conf(X,N)$ are linearly independent by 
Proposition~\ref{prop:FreeBasis}, hence we obtain the following

\begin{corollary}[\cite{Ko2020}]\label{cor:FreeBasis_module}
The free associative conformal algebra $\Conf(X,N)$ is isomorphic 
to $M(X,N)$ as an $A(X)$-module.
\end{corollary}

It follows from the definition of the action of $A(X)$ on $\Conf(X,N)$ that 
every conformal ideal of $\Conf(X,N)$ is an $A(X)$-submodule and vice versa. 
Hence we may replace the study of conformal ideals with the study of 
``ordinary'' submodules.

\begin{example}[\cite{BFK2000}]
 Let us determine the structure of  an associative conformal algebra $C$ 
generated by
 the set $X=\{a\}$ relative to $N=N(a,a)=2$ with one defining relation 
 $a\oo{1} a - \partial (a\oo 0 a)$.
\end{example}
 
The algebra $A(X)$ is generated by $L_n=L_n^a$, $R_n=R_n^a$, and $\partial $
satisfying \eqref{eq:rel_A(X)}. Namely, consider these relations as rewriting 
rules 
\[
L_n \partial \to \partial L_n + n L_{n-1}, \quad
R_n \partial \to \partial R_n + n R_{n-1}, \quad
R_n L_m \to  L_m R_n.
\]
Similarly, define the free conformal algebra $\Conf(X,N)$ as a module over 
$A(X)$ generated by a single element $a$ relative to the following rewriting 
rules \eqref{eq:mod_A(X)}:
\[
\begin{gathered}
L_n a\to 0,\ R_n a  \to 0,\ n\ge 2,
R_1a\to -L_1a,\quad R_0a\to L_0a-\partial L_1a.
\end{gathered}
\]
The compositions \eqref{eq:mod_A(X)_locality}  of these relations include 
\[
L_3L_1\to 0,\ L_3L_0a\to 0,\quad L_2L_1\to 0, \quad L_2L_0\to 2L_1L_1.
\]

The defining relation $a\oo{1} a - \partial (a\oo 0 a)$ is naturally written as
\begin{equation}\label{eq:Exm_BFK}
L_1 a \to \partial L_0 a.
\end{equation}

Consider the composition of $R_2L_1\to L_1R_2$ and \eqref{eq:Exm_BFK} relative 
to $w=R_2L_1a$. 
On the one hand, 
\[
R_2L_1a \to L_1R_2a \to 0,
\]
on the other hand,
\[
R_2L_1a \to R_2\partial L_0 a \to \partial R_2L_0a + 2R_1L_0a \to 2 L_0R_1a \to 
-2L_0L_1a.
\]
Hence, we should add a new rewriting rule 
\begin{equation}\label{eq:BFK-1}
L_0L_1a \to 0.
\end{equation}
The latter has a composition with \eqref{eq:Exm_BFK} relative to $w=L_0L_1a$:
\[
L_0L_1a \to L_0\partial L_0a \to \partial L_0L_0a.
\]
Hence, we should add 
\[
\partial L_0L_0 a\to 0.
\]

Next, consider the composition of $R_1L_1\to L_1R_1$ and \eqref{eq:Exm_BFK} 
relative to $w=R_1L_1a$. 
In a similar way, we obtain that $-L_1L_1a$ and $L_0L_0a$ are connected by a 
(non-oriented) path, 
so add 
\begin{equation}\label{eq:BFK-2.5}
L_1L_1a \to -L_0L_0a
\end{equation}
(the choice of the principal part is voluntary since we have not fixed an order 
on the words).

Let us calculate the composition of $R_0L_1a\to L_1R_0a$ and \eqref{eq:Exm_BFK} 
relative to 
$w=R_0L_1a$ in more details:
\[
\begin{aligned}
R_0L_1a \to L_1R_0a &{} \to L_1L_0a - L_1\partial L_1 a \to L_1L_0-\partial 
L_1L_1a - L_0L_1a \\
 &{}\to L_1L_0a + \partial L_0L_0a \to L_1L_0a, \\
R_0L_1a&{}\to R_0\partial L_0a \to \partial L_0R_0a \to \partial L_0L_0a 
-\partial L_0\partial L_1a \\
 &{}\to  \partial L_0L_0a -\partial^2 L_0 L_1a  \to \partial L_0L_0a\to 0.
\end{aligned}
\]
Hence, we should add
\begin{equation}\label{eq:BFK-2}
L_1L_0a\to 0.
\end{equation}

There exist compositions between \eqref{eq:mod_A(X)_locality} and 
\eqref{eq:Exm_BFK}. 
For example, the composition of $L_2L_1a\to 0 $ and \eqref{eq:Exm_BFK} is 
trivial:
\[
L_2L_1a \to L_2\partial L_0 a \to \partial L_2L_0a + 2L_1L_0a \to 2\partial 
L_1L_1a \to 
-2\partial L_0L_0 a \to 0.
\]
However, the composition of \eqref{eq:Exm_BFK} and $L_3L_1a\to 0$ is not 
trivial:
\[
L_3L_1a \to L_3\partial L_0 a \to \partial L_3L_0a + 3L_2L_0a \to 6L_1L_1a.
\]
Hence, we should add 
\begin{equation}\label{eq:BFK-3}
L_1L_1a\to 0
\end{equation}
and \eqref{eq:BFK-2.5} implies 
\begin{equation}\label{eq:BFK-4}
L_0L_0a\to 0.
\end{equation}

The relations 
\eqref{eq:BFK-1}, \eqref{eq:BFK-2}, \eqref{eq:BFK-3}, \eqref{eq:BFK-4}
along with \eqref{eq:mod_A(X)} and \eqref{eq:mod_A(X)_locality}
form a Gr\"obner--Shirshov basis of $C$ as of $A(X)$-module: all other 
compositions are trivial by 
homogeneity reasons. As a result, the basis of $C$ as of a module over $H=\Bbbk 
[\partial ]$ consists 
of two elements: $a$ and $a\oo 0 a$, all words of degree $\ge 3$ are zero.

\subsection{Universal associative conformal envelopes of Lie conformal 
algebras}\label{subsec_env}

Suppose $L$ is a Lie conformal algebra generated by a set $X$. 
Thus $L$ is a quotient of an appropriate free Lie conformal algebra 
by the ideal generated by a set $\Sigma $ of defining relations stated in terms 
of Lie conformal operations 
$[x\oo{n} y]$. The structure of free Lie conformal algebras was described in 
\cite{Ro1999}.

For a given function $N: X\times X\to \mathbb Z_+$, the universal enveloping 
associative conformal algebra $U(L; X,N)$ of $L$ relative to the locality level 
$N$ on $X$ is defined 
as the quotient of $\Conf (X,N)$ relative to the same defining relations $\Sigma 
$
rewritten by the rules
\[
[x\oo{n} y] = (x\oo{n} y) -\sum\limits_{s\ge 
0}\frac{(-1)^{n+s}}{s!}\partial^{s} (y\oo{n+s} x), 
\] 
where the upper limit of the summation is determined by the Dong Lemma.

The main purpose of this paper is to study universal enveloping 
associative conformal algebras for Kac--Moody conformal algebras. 
The latter are central extensions of current Lie conformal algebras. 
For this particular class of problems, the Gr\"obner--Shirshov bases method 
described above may be slightly modified. 
The main advantage of the modification is that the relations 
\eqref{eq:mod_A(X)_locality}
become not necessary.

Suppose $L$ is a Lie conformal algebra with an $H$-torsion $L_0$ such that the 
torsion-free 
$L_1= L/L_0$ is a free $H$-module (for example, every finite Lie conformal 
algebra has that property). 
Assume $X=X_1\cup X_0$, where $X_1$ is an $H$-basis of $L_1$ and $X_0$ is a 
$\Bbbk $-basis of $X_0$. 
Then the structure of $L$ is completely determined by relations 
\[
f_e(\partial )e = 0,\quad e\in X_0, 
\]
and 
\[
[x\oo{n} y ] = \sum\limits_{z\in X_1} f_{x,y}^{n,z} (\partial )z + 
\sum\limits_{e\in X_0}g_{x,y}^{n,e}(\partial )e,\quad x,y\in X_1,
\]
for appropriate $f_e,f_{x,y}^{n,z}, g_{x,y}^{n,e} \in \Bbbk [\partial ]$.
These relations describe the structure of $L_0$ as of a torsion $H$-module, the 
multiplication table 
in the Lie conformal algebra $L_1$, and the structure of the extension 
\[
0\to L_0 \to L\to L_1\to 0.
\]

Then, for a given function $N:X_1\times X_1\to \mathbb Z_+$, 
the conformal algebra $U(L; X,N)$ may be considered as an ordinary left module 
over the associative algebra 
$A(X;L)$ generated by 
\[
\{\partial, L_n^x, R_n^x \mid n\in \mathbb Z_+, x\in X_1\}
\]
relative to defining relations \eqref{eq:rel_A(X)} (for $a\in X_1$) along with 
the following ones:
\begin{equation}\label{eq:A(X;L)}
L_n^xL_m^y - L_m^yL_n^x - \sum\limits_{s\ge 0}\sum\limits_{z\in X_1} 
L_{n+m-s}^{f_{x,y}^{n,z} (\partial )z}
\end{equation}
where 
$L_n^{\partial z}$ is naturally understood as $-nL_{n-1}^z$. 
The relations \eqref{eq:A(X;L)} reflect the property \eqref{eq:Conf_JacComm} of 
associative conformal algebras. 
So $U(L;X,N)$ is a left module over $A(X;L)$ generated by the entire set $X$ 
relative to the relations
\eqref{eq:mod_A(X)} (for $a,b\in X_1$) together with 
\begin{gather}\label{eq:mod_U(L)}
f_e(\partial ) e,\quad e\in X_0, \\
L_n^a e,\ R_n^a e \quad a\in X_1, \ e\in X_0,\ n\in \mathbb Z_+, \\
R_n^a b - L_n^a b + [a\oo{n} b],\quad a,b \in X_1. 
\end{gather}

Since the defining relations of $A(X;L)$ already form a Gr\"obner--Shirshov 
basis, in order to determine the structure 
of $U(L;X,N)$ one needs to find a confluent system of rewriting rules in this 
$A(X;L)$-module. In the next section, we solve this problem for a Kac--Moody 
conformal algebra.

\section{The Poincar\'e--Birkhoff--Witt Theorem for Kac--Moody conformal 
algebras at $N=3$}\label{sec:PBW-KM}

Let $\mathfrak g$ be a Lie algebra and let $\<\cdot|\cdot \>$ be a bilinear 
symmetric invariant form on $\mathfrak g$ (e.g., the Killing form). 
Then 
$K(\mathfrak g) = (\Bbbk[\partial ]\otimes \mathfrak g) \oplus \Bbbk e$, where 
$\partial e=0$, 
equipped with 
\[
[a\oo\lambda b] = [a,b] + \lambda \<a|b\>, \quad [a\oo\lambda e]=[e\oo\lambda e] 
= 0
\]
for every $a,b\in \mathfrak g$
is a Lie conformal algebra with 1-dimensional torsion $\Bbbk e$
and the torsion-free image isomorphic to $\Cur \mathfrak g$.

Let us fix a linear basis $X_1$ of $\mathfrak g$. Then $X=X_1\cup \{e\}$ is a 
generating set of $K(\mathfrak g)$. The purpose of this section is to calculate 
the Gr\"obner--Shirshov basis for $U = U(K(\mathfrak g); X, N)$ for $N=3$ and 
prove the Poincar\'e--Birkhoff--Witt Theorem for this universal enveloping 
associative conformal algebra.

According to the scheme described in the previous section, $U$ is a module over 
the associative algebra $A=A(X;K(\mathfrak g))$
generated by the set 
$B = \{ \partial, L_n^a, R_n^a\mid a\in X_1,n\in \mathbb Z_+\}$ modulo the 
relations 
\begin{gather}
    L_n^a\partial - \partial L_n^a - nL_{n-1}^a, \quad 
    R_n^a\partial - \partial R_n^a - nR_{n-1}^a, \quad 
    R_n^aL_m^b-L_m^bR_n^a,          \label{eq:U-A-relations1}\\
    L_n^aL_m^b - L_m^bL_n^a - L_{n+m}^{[a,b]}. \label{eq:U-A-relations2}
\end{gather}
The set of generators of $U$ as of an $A$-module is $X=X_1\cup \{e\}$,
and the defining relations of this module are 
\begin{gather}
L_n^ab, R_n^a b\quad n\ge N=3, 
       \label{eq:U-Mod-relations1}  \\
L_n^ae, R_n^ae, \quad n\ge 0, 
            \label{eq:U-Mod-relations2}\\
       R_2^a b - L_2^ba, \quad R_1^ab + L_1^ba -\partial L_2^ba, 
       \label{eq:U-Mod-relations3}\\
R_0^ab - L_0^ba+\partial L_1^b a - \frac{1}{2}\partial^2 L_2^ba,
\label{eq:U-Mod-relations4}\\
R_0^ab-L_0^ab + [a,b], \quad R_1^ab-L_1^ab +\<a|b\>e, \quad 
  R_2^ab-L_2^ab,
   \label{eq:U-Mod-relations5}
\end{gather}
for all $a,b\in X_1$.

In order to translate these defining relations into rewriting rules we need to 
choose a principal monomial in each relation.
The choice of principal parts affects on the resulting system of rewriting rules 
obtained in a process of adding compositions
similar to Example \ref{eq:Exm_BFK}.

We will always choose a principal term in a rewriting rule 
as a leading monomial relative to an appropriate order $\le $ on the monomials 
in the free $\Bbbk\<B\>$-module generated by $X$. 
Namely, suppose the set $X_1$ is linearly well ordered and $e<X_1$.
Induce an order on $B$ by the rule
\[
L_0^a < L_1^a<\partial <L_2^a <\dots <R_0^a <R_1^a < \dots,
\]
assuming $L_n^a<R_m^b$ and $L_n^a<L_n^b$ iff $a<b$, for $a,b\in X_1$
(this ordering turns to be the most convenient for our purpose). 
Extend the order on the set of monomials in $B^*$ by the deg-lex 
principle, i.e., first compare the lengths and then lexicographically.

For two monomials  $ux$ and $vy$ in 
$\Bbbk\<B\> X$, $u,v\in B^*$,
$x,y\in X$, set 
$ux<vy$ iff $(u,x)$ is lexicographically less than $(v,y)$.

Then, eliminating the monomials $R_n^ab$ ($n=0,1,2$)
in \eqref{eq:U-Mod-relations3} and \eqref{eq:U-Mod-relations4}
by means of \eqref{eq:U-Mod-relations5}, we obtain the following 
set of rewriting rules defining $U$ (along with the confluent set of rewriting 
rules for $A$):
\begin{equation}\label{eq:U3-rulesA}
 \begin{gathered}
 \partial L_0^a \to L_0^a\partial ,
 \quad   \partial L_1^a\to L_1^a\partial -L_0^a, \\
 L_n^a\partial \to \partial L_0^a + nL_{n-1}^a,\quad n\ge 2, \\
 R_n^a\partial \to \partial R_n^a + nR_{n-1}^a,\quad n\ge 0, \\
 R_m^aL_n^b \to L_n^bR_m^a,\quad n,m\ge 0, \\
 L_n^aL_m^b \to L_m^bL_n^a + L_{n+m}^{[a,b]},\quad (n,a)>_{lex}(m,b);
 \end{gathered}
\end{equation}
\begin{equation}\label{eq:U3-rulesB}
 \begin{gathered}
L_n^ab\to 0,\quad R_n^a\to 0,\quad n\ge 3, \\
R_2^ab \to L_2^ab,\quad R_1^ab\to L_1^{a}b - \<a|b\>e, \\
R_0^ab \to L_0^ab - [a,b];
 \end{gathered}
\end{equation}
\begin{equation}\label{eq:U3-rulesC}
 \begin{gathered}
 L_2^ab \to L_2^ba,\quad a>b, \\
 \partial L_2^ab \to L_1^ab +L_1^ba -\<a|b\>e, \\
 L_1^a\partial b \to L_1^b\partial a + 3L_0^ab-3L_0^ba -2[a,b],
 \quad a>b.
 \end{gathered}
\end{equation}

\begin{theorem}\label{thm:GSB}
The set of rewriting rules
\eqref{eq:U3-rulesA}--\eqref{eq:U3-rulesC}
along with 
\begin{equation}\label{eq:GSB-Ds}
\begin{gathered}
L_1^a\partial^sb \to L_1^b\partial^sa - (s+2)L_0^b\partial^{s-1} a + 
(s+2)L_0^a\partial^{s-1}b - 2\partial^{s-1}[a, b], \\
\quad s\ge2,\ a>b,
\end{gathered}
\end{equation}
\begin{gather}
L_2^aL_2^bc \to 0, \quad a,b,c\in X_1           \label{eq:GSB-2.2} 
\\
L_1^aL_2^bc \to L_1^bL_2^ca ,\quad b\le c<a,         \label{eq:GSB-1.2}
\\
L_1^aL_2^bc \to L_1^bL_2^ac ,\quad b<a\le c,          \label{eq:GSB-1.2'}
\\
L_1^aL_1^bc \to  L_1^a L_1^c b +  L_0^b L_2^a c - L_0^c L_2^a b + 
L_2^a [c, b] + L_2^b [c,a] + L_2^c [a,b],\quad a\le c<b,    \label{eq:GSB-1.1}
\\
L_1^aL_1^bc \to L_1^c L_1^a b +  L_0^b L_2^c a -  L_0^c L_2^a b + 
L_2^c [a,b] + L_2^a [c, b],\quad c<a\le b,         \label{eq:GSB-1.1'}
\end{gather}
\begin{multline}
L_0^aL_1^bc \to L_0^a L_1^c b +  L_0^b L_1^a c + L_0^c L_1^b a -
L_0^b L_1^c a - L_0^c L_1^a b + L_1^{[c,a]} b +\\+ L_1^{[a,b]} c + 
L_1^{[b,c]} a  - L_1^c [a,b] - L_1^a [b,c] - L_1^b [c,a] + 
\<a|[b, c]\> e,\quad c<b<a,         \label{eq:GSB-0.1}
\end{multline}
is a Gr\"obner--Shirshov basis of 
the universal associative envelope 
$U=U(K(\mathfrak g); X,3)$.
\end{theorem}

\begin{proof}
First, we will show how to derive the rules 
\eqref{eq:GSB-Ds}--\eqref{eq:GSB-0.1} 
as compositions of the initial relations.
Next, we will check the triviality of compositions obtained in further 
iterations. 

Since 
the calculations are routine, we will state 
them in details for several particular cases, 
other cases  are essentially the same and may 
be processed in a similar way.

For the purpose of clarity, we will use a brief notation to 
point a rule applied for rewriting (e.g., $(RL)$ stands 
for $R_m^aL_n^b \to L_n^bR_m^a,n,m\ge 0$, $(\partial L_2)$ 
for $\partial L_2^ab \to L_1^ab +L_1^ba -\<a|b\>e$, etc).

The rule \eqref{eq:GSB-Ds} for $s=2$ 
appears from the intersection of 
$(\partial L_1)$ 
and $(L_1 \partial)$. Then, by induction, the intersection with $(\partial L_1)$ 
produces 
\eqref{eq:GSB-Ds} for  $s> 2$:
\begin{multline*}
\partial L_1^a\partial^s b 
\to 
\partial \big (
L_1^a\partial ^sb - (s+2)L_0^b\partial^{s-1}a + (s+2)L_0^a\partial^{s-1}b 
-2\partial^{s-1}[a,b]
\big) \\
\to
L_1^a\partial ^{s+1}b - (s+3)L_0^b\partial^{s}a + (s+2)L_0^a\partial^{s}b 
-2\partial^{s}[a,b];
\end{multline*}
\begin{equation*}
\partial L_1^a\partial^s b \to 
L_1^a\partial^{s+1} b - L_0^a \partial^s b.
\end{equation*}

The rule \eqref{eq:GSB-2.2} fairly simply derives from $(\partial L_2)$ in 
\eqref{eq:U3-rulesC} and $(L_n \partial)$ for $n = 2$.

The next example of an intersection of $(L_1 \partial)$ and $(RL)$ produces the rest of the required rules. 
On the one hand, we have  
\[
R_n^a L_1^b \partial c 
\overset{(RL)}\to
L_1^b R_n^a \partial c 
\overset{(R\partial)}\to 
L_1^b \partial R_n^a c + n L_1^b R_{n-1}^a c,
\]
on the other hand,
\begin{multline} \nonumber
R_n^a L_1^b \partial c 
\overset{(L_1 \partial)}\to 
R_n^a L_1^c \partial b - 3 R_n^a L_0^c b + 
3 R_n^a L_0^b c - 2 R_n^a [b, c] 
\overset{(RL)}\to 
L_1^c R_n^a \partial b - 3 L_0^c R_n^a b \\
+  3 L_0^b R_n^a c - 2 R_n^a [b, c] 
\overset{(R \partial)}\to 
L_1^c \partial R_n^a b + n L_1^c R_{n-1}^a b - 
3 L_0^c R_n^a b + 3 L_0^b R_n^a c - 2 R_n^a [b, c].
\end{multline}
In order to apply $(L_1\partial)$ we have to assume $b>c$. However, the 
composition obtained by subtracting the right-hand sides of the two expressions 
above is
\begin{equation}\label{eq:L_1d-RL-n}
L_1^b \partial R_n^a c + n L_1^b R_{n-1}^a c 
-
L_1^c \partial R_n^a b - n L_1^c R_{n-1}^a b + 
3 L_0^c R_n^a b - 3 L_0^b R_n^a c + 2 R_n^a [b, c],
\end{equation}
it is (skew-)symmetric relative to the permutation of $b$ and $c$.
Hence, we may assume the relation 
\eqref{eq:L_1d-RL-n} holds on $U$ for every $a,b,c\in X_1$.

For $n \ge 4$ the composition is trivial due to 
the locality. 
For $n = 3$, apply the rules $(R_3)$ and $(R_2)$ to get
 the following:
\begin{equation} \label{eq:Id-1.2}
L_1^b L_{2}^a c - L_1^c L_{2}^a b,
\quad 
a, b, c \in X_1.
\end{equation} 
For a fixed order on $a,b,c\in X_1$, use $(L_2)$ if necessary to obtain 
\eqref{eq:GSB-1.2} or \eqref{eq:GSB-1.2'}. 

Consider \eqref{eq:L_1d-RL-n} for $n=2$. 
For convenience of the exposition, let us split the polynomial into two summands and 
process the summands separately:
\begin{equation}\label{eq:Ld-cap-RL}
\big( L_1^b \partial R_2^a c + 2 L_1^b R_{1}^a c\big )
-\big( L_1^c \partial R_2^a b + 
2 L_1^c R_{1}^a b - 3 L_0^c R_2^a b + 3 L_0^b R_2^a c - 2 R_2^a [b, c] \big ),
\end{equation}
\begin{multline}\nonumber
L_1^b \partial R_2^a c + 2 L_1^b R_{1}^a c 
\overset{(R_1), (R_2)}\to  
L_1^b \partial L_2^a c + 2 L_1^b L_{1}^a c - 2\<a|c\> L_1^b  e
\\
\overset{(\partial L_2), (Le)}\to 
L_1^b L_1^a c + L_1^b L_1^c a - \<c|a\>L_1^b  e 
+ 2 L_1^b L_1^a c 
\overset{(LL), (Le)}\to
3 L_1^b L_{1}^a c + L_1^c L_1^b a + 
L_2^{[b, c]} a, 
\end{multline}
\begin{multline}\nonumber
L_1^c \partial R_2^a b + 2 L_1^c R_{1}^a b - 3 L_0^c R_2^a b + 3 L_0^b R_2^a c 
- 2 R_2^a [b, c] 
\overset{(R_1), (R_2)}\to 
L_1^c \partial L_2^a b + 2 L_1^c L_{1}^a b 
\\
- 2\<a|b\> L_1^c  e - 3 L_0^c L_2^a b + 3 L_0^b L_2^a c - 2 L_2^a [b, c] 
\overset{(\partial L_2), (Le)}\to  
L_1^c L_1^a b + L_1^c L_1^b a - \<b|a\> L_1^c  e 
+ 2 L_1^c L_1^a b \\
- 3 L_0^c L_2^a b + 3 L_0^b L_2^a c - 2 L_2^a [b, c] 
\overset{(Le)}\to 3 L_1^c L_{1}^a b + L_1^c L_1^b a - 3 L_0^c L_2^a b + 
3 L_0^b L_2^a c - 2 L_2^a [b, c]. 
\end{multline}
Therefore,  \eqref{eq:Ld-cap-RL} modulo $(L_2)$ (i.e., $L_2^{[b,c]}a- 
L_2^a[b,c]$) implies 
the following relation: 
\begin{equation} \label{eq:Id-1.1}
- L_1^b L_{1}^a c +  L_1^c L_{1}^a b - L_0^c L_2^a b + 
L_0^b L_2^a c + L_2^a [c, b]
\end{equation} 
for all  $a, b, c \in X_1$.

Now we can switch $a$ and $b$ in \eqref{eq:Id-1.1} and subtract the relation 
obtained from \eqref{eq:Id-1.1}. In this way, we actually apply the rule $(LL)$ 
without fixing an order on $a,b\in X_1$. As a result, we obtain
\begin{multline*}
 - L_1^b L_{1}^a c +  L_1^c L_{1}^a b 
 - L_0^c L_2^a b + 
L_0^b L_2^a c + L_2^a [c, b] \\
+ L_1^a L_{1}^b c - L_1^c L_{1}^b a 
+ L_0^c L_2^b a - L_0^a L_2^b c - L_2^b [c, a] 
\end{multline*}
Let us apply $(LL)$ and $(L_2)$ to write the last relation in 
a more convenient form:
\begin{equation}
 \label{eq:Id-1.1'}
 L_1^c L_{1}^a b - L_1^c L_{1}^b a + 
L_0^b L_2^a c - L_0^a L_2^b c 
+ L_2^c [a, b] + L_2^a [c, b] + L_2^b [a, c].
\end{equation}
Given a fixed order on $a,b,c\in X_1$, use $(LL)$ and $(L_2)$ if necessary to 
obtain \eqref{eq:GSB-1.1} or \eqref{eq:GSB-1.1'}. 

For $n=1$, proceed with \eqref{eq:L_1d-RL-n} in a similar way:
\[
\big (L_1^b \partial R_1^a c + L_1^b R_{0}^a c \big ) -\big ( L_1^c \partial R_1^a b + 
L_1^c R_{0}^a b - 3 L_0^c R_1^a b + 3 L_0^b R_1^a c - 2 R_1^a [b, c]\big ).
\]
On the one summand, we have
\begin{multline}\nonumber
L_1^b \partial R_1^a c + L_1^b R_{0}^a c 
\overset{(R_0), (R_1)}\to
L_1^b \partial L_1^a c + L_1^b L_{0}^a c - L_1^b [a, c]
\overset{(\partial L_1)}\to
L_1^b L_1^a \partial c - L_1^b L_0^a c
\\
+ L_1^b L_{0}^a c - L_1^b [a, c] 
\overset{(L_1 \partial)}\to
L_1^b L_1^c \partial a - 3 L_1^b L_0^c a + 3 L_1^b L_0^a c 
- 2 L_1^b [a, c] - L_1^b [a, c] 
\\ 
\overset{(LL)}\to 
L_1^c L_1^b \partial a + L_2^{[b, c]} \partial a - 3 L_0^c L_1^b a 
- 3 L_1^{[b, c]} a + 3 L_0^a L_1^b c + 3 L_1^{[b, a]} c - 3 L_1^b [a, c],
\end{multline}
whereas on the other
\begin{multline}\nonumber 
L_1^c \partial R_1^a b + L_1^c R_{0}^a b - 3 L_0^c R_1^a b + 
3 L_0^b R_1^a c - 2 R_1^a [b, c] 
\overset{(R_0), (R_1)}\to
L_1^c \partial L_1^a b 
+ L_1^c L_{0}^a b - L_1^c [a, b] 
\\
- 3 L_0^c L_1^a b + 3 L_0^b L_1^a c - 2 L_1^a [b, c] + 2 \<a|[b, c]\> e 
\overset{(\partial L_1)}\to 
L_1^c L_1^a \partial b - L_1^c L_0^a b + L_1^c L_{0}^a b - L_1^c [a, b] 
\\
- 3 L_0^c L_1^a b + 3 L_0^b L_1^a c - 2 L_1^a [b, c] + 2 \<a|[b, c]\> e 
\overset{(L_1 \partial)}\to 
L_1^c L_1^b \partial a - 3 L_1^c L_0^b a + 3 L_1^c L_0^a b - 2 L_1^c [a, b] 
\\
- L_1^c [a, b] - 3 L_0^c L_1^a b + 3 L_0^b L_1^a c - 2 L_1^a [b, c] + 2 \<a|[b, 
c]\> e 
\overset{(LL)}\to 
L_1^c L_1^b \partial a - 3 L_0^b L_1^c a - 3 L_1^{[c, b]} a 
\\
+ 3 L_0^a L_1^c b + 3 L_1^{[c, a]} b - 3 L_1^c [a, b] - 3 L_0^c L_1^a b 
+ 3 L_0^b L_1^a c - 2 L_1^a [b, c] + 2 \<a|[b, c]\> e. 
\end{multline}
Once the expressions subtracted, the terms $L_1^cL_1^b\partial a$ 
cancel among other similar terms, and the result may be rewritten 
via $(L_2 \partial)$ and $(\partial L_2)$.
The resulting expression is symmetric
in $a,b,c\in X_1$, so as we fix the order $c<b<a$  
the principal part of the relation obtained is $L_0^aL_1^bc$, and the rule 
\eqref{eq:GSB-0.1} follows. 

For $n=0$, the relation \eqref{eq:L_1d-RL-n} turns into
\[
L_1^b \partial R_0^a c -\big( L_1^c \partial R_0^a b - 3 L_0^c R_0^a b + 
3 L_0^b R_0^a c - 2 R_0^a [b, c] \big ).
\]
Without loss of generality, assume $b>c$. Then
the two summands in the above relation may be rewritten as follows:
\begin{multline} \nonumber
L_1^b \partial R_0^a c 
\overset{(R_0)}\to 
L_1^b \partial L_0^a c - L_1^b \partial [a, c] 
\overset{(\partial L_0), (LL)}\to 
L_0^a L_1^b \partial c + L_1^{[b, a]} \partial c - L_1^b \partial [a, c] 
\\
\overset{(L_1 \partial)}\to
L_0^a L_1^c \partial b - 3 L_0^a L_0^c b + 3 L_0^a L_0^b c - 2 L_0^a [b, c] +
L_1^c \partial [b, a] - 3 L_0^c [b, a] + 3 L_0^{[b, a]} c -
\\
2 [[b, a], c] - L_1^{[a, c]} \partial b + 3 L_0^{[a, c]} b -
3 L_0^b [a, c] + 2 [b, [a, c]],
\end{multline}
\begin{multline}\nonumber
L_1^c \partial R_0^a b - 3 L_0^c R_0^a b + 3 L_0^b R_0^a c - 2 R_0^a [b, c] 
\overset{(R_0)}\to
L_1^c \partial L_0^a b - L_1^c \partial [a, b] - 3 L_0^c L_0^a b +
\\
3 L_0^c [a, b] + 3 L_0^b L_0^a c - 3 L_0^b [a, c] - 2 L_0^a [b, c] + 2 [a, [b, 
c]] 
\\
\overset{(\partial L_0), (LL)}\to 
L_0^a L_1^c \partial b + L_1^{[c, a]} \partial b -
L_1^c \partial [a, b] - 3 L_0^a L_0^c b - 3 L_0^{[c, a]} b + 3 L_0^c [a, b] 
\\
+
3 L_0^a L_0^b c + 3 L_0^{[b, a]} c - 3 L_0^b [a, c] - 2 L_0^a [b, c] + 2 [a, [b, 
c]]  .
\end{multline}
Hence, the composition \eqref{eq:L_1d-RL-n} for 
$n=0$ is trivial due to the skew-symmetry and the Jacobi identity on 
$[\cdot,\cdot]$.

In order to finish the proof one needs to check that the family of rewriting 
rules obtained is 
complete, i.e., forms a Gr\"obner--Shirshov basis.
Let us consider several intersections as examples, 
other possible intersections may be processed in a similar way.

As a first example, consider the intersection of 
\eqref{eq:GSB-1.2} and $(\partial L_1)$:
\[
\partial L_1^aL_2^bc  
\to 
\partial L_1^bL_2^c a 
\to 
L_1^b\partial L_2^ca -L_0^bL_2^ca  
\to 
L_1^bL_1^ca+L_1^bL_1^ac -L_0^bL_2^c a, 
\]
\begin{multline*}
\partial L_1^aL_2^bc 
\to 
L_1^aL_2^bc-L_0^aL_2^b c 
\to 
L_1^aL_1^bc + L_1^aL_1^cb -L_0^aL_2^bc \\
\to 
L_1^bL_1^ac + L_1^cL_1^ab -L_0^aL_2^bc + L_2^{[a,b]}c + L_2^{[a,c]} b.
\end{multline*}
Subtract the relations obtained to get a composition
\begin{equation}\label{eq:CompExm1}
L_1^cL_1^ab - L_1^bL_1^ca + L_0^bL_2^c a
  -L_0^aL_2^bc + L_2^{[a,b]}c + L_2^{[a,c]} b,
\quad 
b\le c<a.
\end{equation}
If $b<c$ then apply \eqref{eq:GSB-1.1'} to rewrite the composition 
\eqref{eq:CompExm1} into
\[
 L_2^b[c,a] +L_2^c[b,a] + L_2^{[a,b]}c + L_2^{[a,c]} b.
\]
The latter reduces to zero by $(L_2)$.
If $b=c$ in  then apply 
\eqref{eq:GSB-1.1} and $(L_2)$ to reduce \eqref{eq:CompExm1} to zero. 

As a more complicated example,
consider the intersection of
$(RL)$ with \eqref{eq:GSB-1.1}. On the one hand,
\begin{equation}\label{eq:RLL-LHS}
R_n^d L_1^a L_1^b c 
\overset{(RL)}\to L_1^a L_1^b R_n^d c,
\end{equation}
on the other hand,
\begin{multline}\label{eq:RLL-RHS}
R_n^d L_1^a L_1^b c \overset{\eqref{eq:GSB-1.1}}\to
R_n^d L_1^a L_1^c b +  R_n^d L_0^b L_2^a c - R_n^d L_0^c L_2^a b 
+ R_n^d L_2^a [b,c] + R_n^d L_2^b [c,a] + R_n^d L_2^c [a,b] 
\\ 
\overset{(RL)}\to 
L_1^a L_1^c R_n^d b +  L_0^b L_2^a R_n^d c -  L_0^c L_2^a R_n^d b 
+ L_2^a R_n^d [b,c] + L_2^b R_n^d [c,a] + L_2^c R_n^d [a,b].
\end{multline}
Here $a\le c<b$, as in \eqref{eq:GSB-1.1}.

The composition obtained should be considered for 
different $n$'s. 
If $n \ge 3$ then all terms reduce to zero by the locality. 
For $n=2$ we have
\begin{multline}\nonumber
-L_1^a L_1^b R_2^d c + L_1^a L_1^c R_2^d b +  L_0^b L_2^a R_2^d c 
- L_0^c L_2^a R_2^d b + L_2^a R_2^d [c, b] + L_2^b R_2^d [c,a] 
+ L_2^c R_2^d [a,b]
\\
\overset{(R_2)}\to 
- L_1^a L_1^b L_2^d c +  
L_1^a L_1^c L_2^d b +  L_0^b L_2^a L_2^d c -  L_0^c L_2^a L_2^d b 
+ L_2^a L_2^d [c, b] + L_2^b L_2^d [c,a] + L_2^c L_2^d [a,b]\\
\overset{\eqref{eq:GSB-2.2}}\to 
- L_1^a L_1^b L_2^d c + L_1^a L_1^c L_2^d b .
\end{multline}
The latter is trivial modulo \eqref{eq:Id-1.2} and thus reduces to zero by means 
of 
\eqref{eq:GSB-1.2} or \eqref{eq:GSB-1.2'} depending on the order on $b,c,d$.

For $n=1$ we have
\begin{multline}\label{eq:Comp-1_exmp}
-L_1^a L_1^b R_1^d c +
L_1^a L_1^c R_1^d b +  L_0^b L_2^a R_1^d c -  L_0^c L_2^a R_1^d b + L_2^a R_1^d 
[c, b] \\
+ L_2^b R_1^d [c,a] + L_2^c R_1^d [a,b] 
 \overset{(R_1)}\to
-L_1^a L_1^b L_1^d c + L_1^a L_1^c L_1^d b + L_0^b L_2^a L_1^d c \\
-  L_0^c L_2^a L_1^d b 
+ L_2^a L_1^d [c, b] + L_2^b L_1^d [c,a] + L_2^c L_1^d [a,b]
\end{multline}
since the terms containing $e$ annihilate under the action of $L_n$.
Continue reducing \eqref{eq:Comp-1_exmp} with the rules $(LL)$ and 
\eqref{eq:Id-1.1}:
\begin{multline*}
-L_1^a L_0^b L_2^d c + L_1^a L_0^c L_2^d b  
+ L_1^a L_2^d [b, c] +
L_0^b L_1^d L_2^a c 
+ L_0^b L_3^{[a, d]} c - L_0^c L_1^d L_2^a b - L_0^c L_3^{[a, d]} b \\
+ L_1^d L_2^a [c, b] 
+ L_3^{[a, d]} [c, b] + L_1^d L_2^b [c, a] + L_3^{[b, d]} [c, a] + L_1^d L_2^c 
[a,b] 
+ L_3^{[c, d]} [a, b] 
\\
\overset{(LL), (L_3)}\to 
-L_0^b L_1^a L_2^d c - L_1^{[a, b]} L_2^d c 
+ L_0^c L_1^a L_2^d b + L_1^{[a, c]} L_2^d b 
+ L_1^a L_2^d [b, c] 
+ L_0^b L_1^d L_2^a c \\
- L_0^c L_1^d L_2^a b - L_1^d L_2^a [b, c] - L_1^d L_2^b [a, c] + L_1^d L_2^c 
[a,b]
\overset{(L_2), \eqref{eq:Id-1.2}}\to 0.
\end{multline*}
For $n=0$, first continue reducing \eqref{eq:RLL-LHS} and \eqref{eq:RLL-RHS} 
as follows:
\[
L_1^a L_1^b R_0^d c 
\overset{(R_0)}\to 
L_1^a L_1^b L_0^d c - L_1^a L_1^b [d, c]  
\overset{(LL)}\to
L_0^d L_1^a L_1^b c + L_1^a L_1^{[b, d]} c + L_1^{[a, d]} L_1^b c - L_1^a L_1^b 
[d, c],
\]
\begin{multline}\nonumber
L_1^a L_1^c R_0^d b +  L_0^b L_2^a R_0^d c - L_0^c L_2^a R_0^d b 
+ L_2^a R_0^d [c, b] + L_2^b R_0^d [c,a] + L_2^c R_0^d [a,b] 
\overset{(R_0)}\to 
L_1^a L_1^c L_0^d b -
\\
- L_1^a L_1^c [d, b] + L_0^b L_2^a L_0^d c - L_0^b L_2^a [d, c] 
- L_0^c L_2^a L_0^d b + L_0^c L_2^a [d, b] + L_2^a L_0^d [c, b] - L_2^a [d, [c, 
b]] 
\\
+ L_2^b L_0^d [c,a] - L_2^b [d, [c,a]] + L_2^c L_0^d [a,b] - L_2^c [d, [a,b]] 
\overset{(LL)}\to  
L_0^d L_1^a L_1^c b + L_1^a L_1^{[c, d]} b + L_1^{[a, d]} L_1^c b 
\\
- L_1^a L_1^c [d, b] 
+ L_0^d L_0^b L_2^a c + L_0^{[b, d]} L_2^a c + L_0^b L_2^{[a, d]} c - L_0^b 
L_2^a [d, c] - 
L_0^d L_0^c L_2^a b - L_0^{[c, d]} L_2^a b - L_0^c L_2^{[a, d]} b 
\\
+ L_0^c L_2^a [d, b] + L_0^d L_2^a [c, b] + L_2^{[a, d]} [c, b] - L_2^a [d, [c, 
b]] 
+ L_0^d L_2^b [c,a] + L_2^{[b, d]} [c,a] - L_2^b [d, [c,a]] \\
+ L_0^d L_2^c [a,b] 
+ L_2^{[c, d]} [a,b] - L_2^c [d, [a,b]].
\end{multline}
Now subtract the relations obtained, apply the Jacobi identity to rewrite 
\[
 \begin{gathered}
L_2^a[d,[c,b]] = -L_2^a[c,[b,d]] - L_2^a[[c,d],b],
\quad
L_2^b[d,[c,a]] = -L_2^b[c,[a,d]] - L_2^b [[c,d],a],\\
L_2^c [d,[a,b]] = -L_2^c [a,[b,d]] - L_2^c [[a,d],b],
\end{gathered}
\]
and rearrange the summands 
gathering them into four groups: one starting with $L_0^d$ and other three containing 
the Lie brackets $[a,d]$, $[b,d]$, and $[c,d]$:  
\begin{multline} \nonumber
\big (L_0^d L_1^a L_1^b c - L_0^d L_1^a L_1^c b + L_0^d L_0^c L_2^a b - L_0^d 
L_0^b L_2^a c 
+ L_0^d L_2^a [b,c] + L_0^d L_2^b [a,c] + L_0^d L_2^c [b,a]\big )\\
-
\big (L_1^{[a, d]} L_1^c b - L_1^{[a, d]} L_1^b c + L_0^b L_2^{[a, d]} c - L_0^c 
L_2^{[a, d]} b 
+ L_2^{[a, d]} [c, b] + L_2^b [c, [a, d]] + L_2^c [[a, d], b]\big ) \\
-
\big (L_1^a L_1^c [b, d] - L_1^a L_1^{[b, d]} c + L_0^{[b, d]} L_2^a c -  L_0^c 
L_2^a [b, d] 
+ L_2^a [c, [b, d]] + L_2^{[b, d]} [c,a] + L_2^c [a, [b, d]]\big )  \\
-
\big (L_1^a L_1^{[c, d]} b - L_1^a L_1^b [c, d] + L_0^b L_2^a [c, d] - L_0^{[c, 
d]} L_2^a b 
+ L_2^a [[c, d], b] + L_2^b [[c, d], a] + L_2^{[c, d]} [a,b]\big ).
\end{multline}
The first group is exactly \eqref{eq:Id-1.1} under the action of $L_0^d$, other 
three groups coincide with \eqref{eq:Id-1.1'}. Therefore, the composition 
reduces to zero by means of \eqref{eq:GSB-1.1} and \eqref{eq:GSB-1.1'}.
\end{proof}

Observe that the principal parts of the rewriting rules from 
the Gr\"obner--Shirshov basis found in 
Theorem \ref{thm:GSB}
do not depend on the multiplication table of the original
Lie algebra~$\mathfrak g$ as well as on the choice of 
the form $\langle\cdot \mid \cdot \rangle $.
In particular, if $\mathfrak g$ is an abelian Lie algebra 
and $\langle x\mid y\rangle =0$ for all $x,y\in \mathfrak g$
then $K(\mathfrak g)$ is an abelian Lie conformal algebra
and $U$ coincides with the 1-dimensional split null extension
\[
0\to \Bbbk e \to U \to \Com\Conf(X_1,N=3) \to 0
\]
of the free  commutative conformal algebra
$\Com\Conf(X_1,N=3)$
generated by a linear basis $X_1$ 
of $\mathfrak g$ relative to 
the locality function $N(x,y)=3$, $x,y\in X_1$.

Hence, the linear basis of $\Com\Conf(X_1,N=3)$ 
consists of all those conformal monomials  
described in Proposition~\ref{prop:FreeBasis} that are 
terminal relative to the rewriting rules stated 
in Theorem \ref{thm:GSB}. 

\begin{corollary}\label{cor:ComConf3}
Let $Y$ be a linearly ordered set. 
The linear basis of the free commutative conformal algebra
$\Com\Conf(Y, N=3)$ 
generated by $Y$ relative to the locality function $N=3$
consists of the following conformal monomials:
\[
\begin{gathered}
L_0^{x_1}\dots L_0^{x_n} L_1^{y_1} \dots L_1^{y_m} L_2^z u,
\quad 
x_i \le x_{i+1},\ y_i \le y_{i+1} \le z \le u,\ 
n,m \ge 0, 
\\
L_0^{x_1}\dots L_0^{x_n} L_1^{y_1} \dots L_1^{y_m} \partial^s z, 
\quad 
x_i \le x_{i+1}, \ y_i \le y_{i+1} \le z,\ n,m \ge 0,
\ s \ge 1, 
\\
L_0^{x_1}\dots L_0^{x_n} L_1^{y_1} \dots L_1^{y_m} z, 
\quad 
x_i \le x_{i+1},\ y_i \le y_{i+1} \le z,\ 
n \ge 0,\ m = 0 \text{ or } m \ge 2, 
\\
L_0^{x_1}\dots L_0^{x_n} L_1^y z, \quad 
x_i \le x_{i+1},\ x_n \le y \text{ or } y \le z,\ n \ge 0, 
\end{gathered}
\]
where $x_i,y_i,z,u\in Y$.

If $X=X_1\cup \{e\}$, $Y=X_1$ then the above monomials 
together with the torsion element $e$  form a linear basis of $U=U(K(\mathfrak 
g); X, N=3)$.
\end{corollary}

The conformal algebra $U$ has a natural filtration by 
degree in~$X$. Denote $\mathrm{gr}\,U$ the corresponding 
associated graded conformal algebra.

Note that every rule in Theorem \eqref{thm:GSB}
has the following property: all terms of highest degree in $X$ do not depend on 
$[\cdot,\cdot]$ and $\langle\cdot \mid \cdot \rangle $.
Therefore, if we choose two basic monomials of $U$ described by Corollary 
\ref{cor:ComConf3} and rewrite their conformal product as a linear combination 
of 
basic monomials then the terms of highest degree in the expression obtained 
would be the same as we get for the product of the same monomials in 
$\Com\Conf(X_1,N=3)$.

As a result, we obtain the following analogue of the Poincar\'e--Birkhoff--Witt 
Theorem.

\begin{corollary}\label{cor:PBW-KM}
For every Lie algebra $\mathfrak g$ and for every bilinear symmetric invariant 
form $\<\cdot,\cdot \>$ on $\mathfrak g$ 
we have 
\[
\mathrm{gr}\,U(K(\mathfrak g); X, N=3)
\simeq \Com\Conf(X_1,N=3)\oplus \Bbbk e.
\]
Here $X=X_1\cup \{e\}$ as above, $X_1$ is a basis of $\mathfrak g$.
In particular,  
\[
\mathrm{gr}\,U(\Cur \mathfrak g; X_1, N=3)
\simeq \Com\Conf(X_1,N=3).
\]
\end{corollary}

Let us derive another corollary of Theorem \ref{thm:GSB}.
Since the case of locality $N = 3$ has already been considered
 we may add
a new series of relations to $U(K(\mathfrak g); X, N=3)$
to get a Gr\"obner--Shirshov basis of the universal enveloping associative 
algebra $U(K(\mathfrak g); X, N=2)$.
Namely, it is enough to add relations
\[
L_2^ab\to 0,\quad R_2^ab\to 0,\quad a,b\in X_1
\]
to the Gr\"obner--Shirshov basis calculated in Theorem \ref{thm:GSB}
and compute all intersections.
Since there are no essentially new manipulations, let us just enlist the 
resulting relations.

\begin{theorem}\label{thm:GSB-2}
The set of rewriting rules \eqref{eq:U3-rulesA} along with 
\begin{gather}
L_n^a b \to 0, \quad n \ge 2,  \label{eq:locL2}
\\
R_n^a b \to 0, \quad n \ge 2,  \label{eq:locR2}
\\
R_1^a b \to L_1^a b - \<a|b\> e,  \label{eq:R1-L1}
\\
R_0^a b \to L_0^a b - [a, b],  \label{eq:R0-L0}
\\
L_1^a b \to - L_1^b a + \<a|b\>e, \quad b < a,  \label{eq:L1-L1}
\\
L_1^a a \to \frac{1}{2} \<a|a\>e,  \label{eq:L1'}
\\
L_1^a \partial^s b \to 2 L_0^a \partial^{s-1} b - L_0^b \partial^{s-1} a -
\partial^{s-1}[a, b],  \quad s\ge 1,\label{eq:L1D}
\\
L_1^a L_1^b c \to 0, \quad a \le b < c,  \label{eq:L1L1}
\\
\begin{gathered}
L_0^a L_1^b c \to L_0^c L_1^b a - L_0^b L_1^c a + L_1^a [c, b] + L_1^c [b, a] + 
L_1^b [a, c] + \<a|[b, c]\>e\\
 b < c < a. 
\end{gathered}\label{eq:L0L1}
\end{gather}
is a Gr\"obner--Shirshov basis of $U=U(K(\mathfrak g); X, 2)$.
\end{theorem}

Here the identities \eqref{eq:locL2}--\eqref{eq:L1'} are slightly adjusted
\eqref{eq:U3-rulesB}, \eqref{eq:U3-rulesC}; and \eqref{eq:L1D}--\eqref{eq:L0L1} 
appear as compositions of intersection.

The list of rewriting rules obtained in Theorem \ref{thm:GSB-2} has the same 
properties as in Theorem~\ref{thm:GSB}: the principal parts as well as all 
summands of maximal degree in each rule do not depend neither on a 
multiplication table on $\mathfrak g$ nor on 
$\<\cdot ,\cdot \>$. Hence, the analogue of Corollary \ref{cor:PBW-KM} holds 
for $N=2$. In order to complete the description of $U(K(\mathfrak g); X, N=2)$ 
we need to state explicitly the set of reduced monomials in the free commutative 
conformal algebra of locality $N=2$.

\begin{corollary}
The following monomials form a linear basis of 
$\Com\Conf(Y,N=2)$, where $Y$ is a linearly ordered set:
\[
\begin{gathered}
L_0^{x_1}\dots L_0^{x_n}\partial ^s z, 
  \quad x_1\le \dots \le x_n,\ s\ge 0, \\
L_0^{x_1}\dots L_0^{x_n} L_1^y z, 
  \quad x_1\le \dots \le x_n\le z, \ y<z,
\end{gathered}
\]
for $x_i,y,z\in Y$.
\end{corollary}

\end{document}